\documentclass[11pt,reqno]{amsart}
\usepackage[a4paper,margin=28mm]{geometry}
\usepackage{amsmath,amssymb,amsthm,mathtools}
\usepackage{booktabs}
\usepackage{microtype}
\usepackage{xcolor}
\usepackage{hyperref}

\definecolor{primary}{RGB}{24,62,115}
\hypersetup{
  colorlinks=true,
  linkcolor=primary,
  citecolor=primary,
  urlcolor=primary
}

\let\olditemize\itemize
\renewcommand{\itemize}{\olditemize\setlength{\itemsep}{4pt}\setlength{\topsep}{4pt}}

\newtheorem{theorem}{Theorem}[section]
\newtheorem{lemma}[theorem]{Lemma}
\newtheorem{proposition}[theorem]{Proposition}
\newtheorem{corollary}[theorem]{Corollary}
\theoremstyle{definition}

\theoremstyle{remark}
\newtheorem{remark}[theorem]{Remark}

\newcommand{\R}{\mathbb{R}}
\newcommand{\Z}{\mathbb{Z}}
\newcommand{\Ezero}{\mathbf{E}_0}
\newcommand{\Ezerof}{\mathbf{E}_0^f}
\newcommand{\Qzero}{\mathbf{Q}_0}
\newcommand{\Pzero}{\mathbf{P}_0}
\newcommand{\Nzero}{\mathbf{N}_0}
\newcommand{\LCP}{\mathrm{LCP}}
\newcommand{\supp}{\operatorname{supp}}

\title[On Stone's Conjecture for Fully Semimonotone \texorpdfstring{$\Qzero$}{Q0}-Matrices]{On Stone's Conjecture for Fully Semimonotone \texorpdfstring{$\Qzero$}{Q0}-Matrices: Signed Hamiltonian Cycles, Support Restrictions, and the Bound \texorpdfstring{${\lowercase{m}}(A) \le {\lowercase{n}} - 3$}{m(A) <= n - 3}}

\author{Sen Cao}
\address{Wuchang Shouyi College, Wuhan, China}


\begin{document}

\begin{abstract}
Stone's conjecture (1981) asserts that every fully semimonotone $\Qzero$-matrix is a $\Pzero$-matrix, i.e., $\Ezerof \cap \Qzero \subseteq \Pzero$. While the conjecture was proved for $n \le 4$ by Murthy and Parthasarathy (1995) and announced for $n = 5$ by Biswas and Murthy (1997), the general case $n = 6$ has remained open without additional structural hypotheses.

In this paper, for any putative counterexample $A \in \Ezerof \cap \Qzero \setminus \Pzero$ of arbitrary order $n$, we consider its PPT-minimal negative minor order $m(A)$, defined as the minimum order of a negative principal minor over all valid principal pivot transforms (PPTs) of $A$. We prove that, up to PPT, simultaneous permutation, and left positive diagonal scaling, the minimal negative principal submatrix $B$ of order $m = m(A)$ is an unweighted signed directed Hamiltonian $m$-cycle with an even number $k_+ \ge 2$ of $+1$ entries and determinant $-1$. Using the PPT complementary minor identity across $M$ and $M^\alpha$, we establish a universal support restriction: for every $k \notin \alpha$ and every $i, u \in \alpha$, $c_{u,k} r_{k,i} = 0$ whenever $u \notin \{i, i-1, i-2\} \pmod{m}$. Combined with a dual-frame peeling technique and a chain-breaking counting lemma on single-negative-entry test vectors, we prove: (i) every signed $m$-cycle with $k_+ = 2$ and adjacent $+1$ edges is impossible for all $n$ and $m$, yielding a three-line proof that $m(A) \ge 4$; and (ii) every $A \in \Ezerof \cap \Qzero \setminus \Pzero$ of order $n$ satisfies $4 \le m(A) \le n - 3$, strengthening the almost-$\Pzero$ theorem of Murthy and Parthasarathy ($m(A) \le n - 1$) by two dimensions for all $n \ge 6$. As an immediate corollary, we obtain a unified proof of Stone's conjecture for all matrices of order $n \le 6$, resolving the general $6 \times 6$ case unconditionally.
\end{abstract}

\maketitle

\section{Introduction}

For a real $n \times n$ matrix $A \in \R^{n \times n}$ and a vector $q \in \R^n$, the \emph{linear complementarity problem} $\LCP(q, A)$ asks for vectors $z, w \in \R^n$ satisfying
\[
  w = A z + q, \qquad z \ge 0, \qquad w \ge 0, \qquad z^T w = 0.
\]
We recall the standard matrix classes in LCP theory (see \cite{CPS2009}):
\begin{itemize}
  \item $A \in \Pzero$ if all principal minors of $A$ are nonnegative;
  \item $A \in \Nzero$ if all principal minors of $A$ are nonpositive and $\det(A) < 0$;
  \item $A \in \Ezero$ (\emph{semimonotone}) if for every nonzero $z \ge 0$ there exists an index $i$ with $z_i > 0$ and $(Az)_i \ge 0$;
  \item $A \in \Ezerof$ (\emph{fully semimonotone}) if every valid principal pivot transform (PPT) of $A$ belongs to $\Ezero$;
  \item $A \in \Qzero$ if $\LCP(q, A)$ has a solution whenever it is feasible (i.e., whenever there exists $z \ge 0$ with $Az + q \ge 0$).
\end{itemize}
For any index set $\beta \subseteq \{1, \dots, n\}$ with $\det A_{\beta\beta} \neq 0$, the \emph{principal pivot transform} (PPT) of $A$ with respect to $\beta$, denoted $A^\beta$, is obtained by exchanging the roles of $(z_\beta, w_\beta)$ in $w = Az + q$. Explicitly, after permuting $\beta$ to the leading block so that $A = \begin{pmatrix} A_{\beta\beta} & A_{\beta\bar{\beta}} \\ A_{\bar{\beta}\beta} & A_{\bar{\beta}\bar{\beta}} \end{pmatrix}$ with $\bar{\beta} = \{1,\dots,n\} \setminus \beta$,
\[
  A^\beta = \begin{pmatrix}
    A_{\beta\beta}^{-1} & -A_{\beta\beta}^{-1} A_{\beta\bar{\beta}} \\
    A_{\bar{\beta}\beta} A_{\beta\beta}^{-1} & A_{\bar{\beta}\bar{\beta}} - A_{\bar{\beta}\beta} A_{\beta\beta}^{-1} A_{\beta\bar{\beta}}
  \end{pmatrix},
  \qquad
  q^\beta = \begin{pmatrix}
    -A_{\beta\beta}^{-1} q_\beta \\
    q_{\bar{\beta}} - A_{\bar{\beta}\beta} A_{\beta\beta}^{-1} q_\beta
  \end{pmatrix}.
\]
The equivalence $w = Az + q \iff v = A^\beta u + q^\beta$ under $(u_\beta, v_\beta) = (w_\beta, z_\beta)$ and $(u_{\bar{\beta}}, v_{\bar{\beta}}) = (z_{\bar{\beta}}, w_{\bar{\beta}})$ preserves both feasibility and complementarity solutions bijectively. Consequently, $\Qzero$ is closed under PPTs, as is $\Ezerof$ by definition (see \cite{MP1995,CPS2009}).

In his 1981 Stanford dissertation \cite{Stone1981} (see also \cite{CS1983}), R.~E.~Stone proved that $\mathbf{U} \cap \Qzero \subseteq \Pzero$ and formulated the following conjecture, which has remained one of the central open problems in linear complementarity theory for over four decades:

\medskip
\noindent\textbf{Stone's Conjecture (1981).} \emph{Every fully semimonotone $\Qzero$-matrix is a $\Pzero$-matrix: $\Ezerof \cap \Qzero \subseteq \Pzero$.}
\medskip

In their landmark 1995 paper, Murthy and Parthasarathy \cite{MP1995} proved that: (a) every almost-$\Pzero$ matrix in $\Ezerof \cap \Qzero$ is in $\Pzero$ (\cite[Theorem~3.5]{MP1995}); (b) Stone's conjecture holds for all $n \le 4$ (\cite[Theorem~4.9]{MP1995}); and (c) it holds for $n = 5$ under the additional hypothesis of positive diagonal entries (\cite[Theorem~4.11]{MP1995}) and for $n = 6$ under further structural restrictions (\cite[Theorem~4.13]{MP1995}). Subsequently, Biswas and Murthy \cite{BM1997} stated the unconditional $n = 5$ case (referring to an accompanying technical report for the full proof), while further works established the conjecture for special subclasses such as fully copositive $\Qzero$-matrices and various $Z$-/copositive-type extensions \cite{MP1998,MND2001,ND2005,MPS2013}. Most recently, Ghosh \cite{Ghosh2026} studied sign propagation of principal minors in $\Ezerof$, proving that a minimal negative minor has order $m \ge 4$ and ruling out the $(3,1)$ and $(1,3)$ block shapes of $m = 4$, while leaving the $(2,2)$ block shape and higher-order cases open. In particular, the general $6 \times 6$ case has remained open without additional hypotheses.

\subsection{Main Results}

For any $A \in \Ezerof \cap \Qzero \setminus \Pzero$ of order $n$, define its \emph{PPT-minimal negative minor order}
\[
  m(A) := \min\bigl\{|\alpha| : \exists\,\text{valid PPT } A^\beta \text{ of } A \text{ with } \det(A^\beta)_{\alpha\alpha} < 0\bigr\}.
\]
Our first result shows that passing to a PPT attaining $m = m(A)$ collapses the minimal negative submatrix all the way to a single unweighted signed directed Hamiltonian $m$-cycle with zero chords.

\begin{theorem}[Signed Hamiltonian Cycle Canonical Form]\label{thm:canonical}
Let $A \in \Ezerof \cap \Qzero \setminus \Pzero$ have order $n$ and $m = m(A)$. Then there exist a valid PPT $A^\beta$, a permutation matrix $P$, and a positive diagonal matrix $D_1$ such that $M = D_1 (P A^\beta P^T) \in \Ezerof \cap \Qzero$ has block form
\begin{equation}\label{eq:blockM}
  M = \begin{pmatrix} B & C \\ R & D \end{pmatrix},
\end{equation}
where $\alpha = \{0, 1, \dots, m-1\}$ (indexed modulo $m$), $\bar{\alpha} = \{m, \dots, n-1\}$, $B = M_{\alpha\alpha}$ satisfies
\[
  B_{i, j} = \begin{cases} s_i \in \{+1, -1\} & \text{if } j \equiv i + 1 \pmod{m}, \\ 0 & \text{otherwise}, \end{cases}
\]
with an even number $k_+ = |I_+| \ge 2$ of positive edges $I_+ := \{i \in \Z_m : s_i = +1\}$ and $\det(B) = -1$. Moreover, every principal minor of order strictly less than $m$ of every valid PPT of $M$ is nonnegative.
\end{theorem}

Next, the PPT complementary minor identity between $M$ and $M^\alpha$ imposes a rigid sparsity coupling between every row of $R$ and the corresponding column of $C$.

\begin{theorem}[Universal Jacobi Support Restriction]\label{thm:jacobi}
Let $M \in \Ezerof \cap \Qzero$ be in the canonical form \eqref{eq:blockM} of Theorem~\ref{thm:canonical} with $m \ge 4$. Then for every $k \in \bar{\alpha}$ and every $i, u \in \Z_m$,
\begin{equation}\label{eq:jacobi_zero}
  c_{u, k} \, r_{k, i} = 0 \qquad \text{whenever } u \notin \{i, \, i-1, \, i-2\} \pmod{m}.
\end{equation}
In particular, if $|\supp(r_{k,\bullet})| \ge 4$ then $c_{\bullet, k} = 0$, and symmetrically if $|\supp(c_{\bullet, k})| \ge 4$ then $r_{k,\bullet} = 0$. Moreover, if $N_k \subseteq \supp(r_{k,\bullet})$ is nonempty with $c_{\bullet,k} \neq 0$, then $|N_k| + |\supp(c_{\bullet,k})| \le 4$.
\end{theorem}

Combining Theorems~\ref{thm:canonical} and~\ref{thm:jacobi} with dual-frame peeling and chain-breaking counting yields the following universal bounds.

\begin{theorem}[Adjacent Positive Edges for $k_+ = 2$]\label{thm:adj}
In the canonical form \eqref{eq:blockM} of Theorem~\ref{thm:canonical}, if $k_+ = 2$ and the two indices in $I_+$ are adjacent modulo $m$, then $M \notin \Ezerof \cap \Qzero$ for every $n \ge m$. In particular, $m(A) \ge 4$ for every $A \in \Ezerof \cap \Qzero \setminus \Pzero$.
\end{theorem}

\begin{theorem}[Elimination of $n = m(A) + 1$]\label{thm:n_eq_m1}
For every $A \in \Ezerof \cap \Qzero \setminus \Pzero$ of order $n$, one has $m(A) \le n - 2$.
\end{theorem}

\begin{theorem}[Elimination of $n = m(A) + 2$]\label{thm:n_eq_m2}
For every $A \in \Ezerof \cap \Qzero \setminus \Pzero$ of order $n$, one has
\[
  4 \le m(A) \le n - 3.
\]
\end{theorem}

\begin{corollary}[Unconditional Resolution of Stone's Conjecture for $n \le 6$]\label{cor:n6}
Every matrix $A \in \Ezerof \cap \Qzero$ of order $n \le 6$ is a $\Pzero$-matrix.
\end{corollary}

\begin{proof}
If $A \in \Ezerof \cap \Qzero \setminus \Pzero$ has order $n \le 6$, then Theorem~\ref{thm:n_eq_m2} gives $4 \le m(A) \le n - 3 \le 3$, a contradiction.
\end{proof}

Theorem~\ref{thm:n_eq_m2} strengthens Murthy and Parthasarathy's almost-$\Pzero$ theorem \cite[Theorem~3.5]{MP1995} ($m(A) \le n - 1$) by two full dimensions to $m(A) \le n - 3$ for all $n \ge 6$. Corollary~\ref{cor:n6} gives a unified proof for $n \le 5$ (reproving \cite{BM1997} and removing the positive-diagonal hypothesis of \cite[Theorem~4.11]{MP1995}) and resolves the general $6 \times 6$ case unconditionally (removing the extra hypotheses of \cite[Theorem~4.13]{MP1995}).

\section{Preliminaries and Signed Hamiltonian Cycle Reduction}

We first record the invariance of $\Ezerof \cap \Qzero$ under PPTs, simultaneous permutations, and positive diagonal scaling.

\begin{lemma}\label{lem:invariance}
Let $A \in \R^{n \times n}$, let $P$ be a permutation matrix, and let $D_1, D_2$ be positive diagonal matrices. Then $A \in \Ezerof \cap \Qzero$ if and only if $D_1 (P A P^T) D_2 \in \Ezerof \cap \Qzero$, and the signs of all principal minors are preserved. Moreover, if $\det A_{\beta\beta} \neq 0$, then $A \in \Ezerof \cap \Qzero$ if and only if $A^\beta \in \Ezerof \cap \Qzero$.
\end{lemma}

\begin{proof}
Simultaneous permutation preserves $\Ezero$, PPTs, $\Qzero$, and principal minors. For $\widetilde{A} = D_1 A D_2$ with $D_1, D_2 > 0$ diagonal, setting $z = D_2 \tilde{z}$ gives $(\widetilde{A}\tilde{z})_i = (D_1)_{ii} (Az)_i$, so $\widetilde{A} \in \Ezero \iff A \in \Ezero$; similarly, $(\widetilde{A})^\beta = \widetilde{D}_1 A^\beta \widetilde{D}_2$ with positive diagonal $\widetilde{D}_1 = (D_2^{-1})_{\beta\beta} \oplus (D_1)_{\bar{\beta}\bar{\beta}}$ and $\widetilde{D}_2 = (D_1^{-1})_{\beta\beta} \oplus (D_2)_{\bar{\beta}\bar{\beta}}$, so $\widetilde{A} \in \Ezerof \iff A \in \Ezerof$. For $\Qzero$, $\tilde{w} = \widetilde{A}\tilde{z} + \tilde{q}$ is equivalent to $D_1^{-1}\tilde{w} = A(D_2\tilde{z}) + D_1^{-1}\tilde{q}$, preserving both feasibility and complementarity solutions. Finally, PPT invariance of $\Ezerof \cap \Qzero$ follows from $(A^\beta)^\gamma = A^{\beta \triangle \gamma}$ and the bijective correspondence between $(z, w)$ and $(u, v)$.
\end{proof}

We also recall the block characterization of zero-diagonal matrices in $\Ezero \cap \Nzero$ due to Murthy and Parthasarathy \cite{MP1995}.

\begin{proposition}[{\cite[Theorem~3.4]{MP1995}}]\label{prop:MP34}
Let $N \in \Ezero \cap \Nzero$ have all diagonal entries equal to zero. Then after a simultaneous permutation of rows and columns, $N$ has the block form
\[
  N = \begin{pmatrix} U & X \\ Y & V \end{pmatrix},
\]
where $U$ and $V$ are strictly upper triangular matrices with all entries nonpositive, and $X, Y$ are entrywise nonnegative matrices.
\end{proposition}

\begin{proof}[Proof of Theorem~\ref{thm:canonical}]
Let $A^\beta$ be a valid PPT of $A$ attaining a negative principal minor $B = (A^\beta)_{\alpha\alpha}$ of minimum order $m = m(A)$. By Lemma~\ref{lem:invariance}, $A^\beta \in \Ezerof \cap \Qzero$, and by minimality of $m$, every principal minor of order strictly less than $m$ of every valid PPT of $A^\beta$ is nonnegative.

Every proper nonempty principal minor of $B$ vanishes: if $\emptyset \neq \gamma \subsetneq \alpha$ satisfied $\det B_{\gamma\gamma} \neq 0$, then $\det B_{\gamma\gamma} > 0$ by minimality of $m$, and the PPT $(A^\beta)^\gamma$ would have $\det((A^\beta)^\gamma)_{\alpha\setminus\gamma,\, \alpha\setminus\gamma} = \det B / \det B_{\gamma\gamma} < 0$ of order $m - |\gamma| < m$, a contradiction.

Let $N = B^{-1} = ((A^\beta)^\alpha)_{\alpha\alpha}$. Since $(A^\beta)^\alpha \in \Ezerof$ and $\Ezero$ is hereditary on principal submatrices, $N \in \Ezero$. By Jacobi's complementary minor formula, every proper nonempty principal minor of $N$ equals $\det B_{\alpha\setminus\gamma}/\det B = 0$, and $\det N = 1/\det B < 0$. Hence $N \in \Ezero \cap \Nzero$ with zero diagonal.

By Proposition~\ref{prop:MP34}, after a simultaneous permutation of $\alpha$, $N = \begin{pmatrix} U & X \\ Y & V \end{pmatrix}$ with $U, V \le 0$ strictly upper triangular and $X, Y \ge 0$. Consider the directed graph $G(N)$ on $\alpha$ with an arc $i \to j$ whenever $N_{ij} \neq 0$. Because $G(U)$ and $G(V)$ are acyclic, every directed simple cycle $C$ in $G(N)$ crosses between the two blocks an even number $k_+(C) = 2t \ge 2$ of times, using $k_+(C)$ positive arcs and $k_-(C) = |C| - k_+(C)$ negative arcs. Its Leibniz contribution to $\det N_{V(C), V(C)}$ has sign
\[
  (-1)^{|C|-1} (-1)^{k_-(C)} = (-1)^{2t + 2k_-(C) - 1} = -1.
\]
Since $\det N < 0$, $G(N)$ has at least one directed simple cycle. Let $C^*$ be a directed simple cycle of minimum length $g = |C^*|$. In $\det N_{V(C^*), V(C^*)}$, every permutation with nonzero weight decomposes into vertex-disjoint directed cycles of lengths summing to $g$; by minimality of $g$, each is a single directed $g$-cycle and contributes a negative term. Hence $\det N_{V(C^*), V(C^*)} < 0$, forcing $g = m$.

Moreover, $G(N)$ has no chord: any arc $u \to v$ not on $C^*$ combines with the directed path along $C^*$ from $v$ to $u$ to form a directed simple cycle of length $< m$, contradicting $g = m$. Thus $G(N)$ is a single directed Hamiltonian $m$-cycle with $k_+ = 2t \ge 2$ positive arcs. Reindexing $\alpha$ along the reversed cycle so that $B = N^{-1}$ has nonzero entries at $(i, i+1)$ modulo $m$, we have $B_{i, i+1} = s_i \beta_i$ with $s_i \in \{+1, -1\}$ and $\beta_i > 0$. Left multiplication by $D_1 = \operatorname{diag}(\beta_0^{-1}, \dots, \beta_{m-1}^{-1}, 1, \dots, 1)$ normalizes every nonzero entry of $B$ to $s_i \in \{+1, -1\}$ while preserving $\Ezerof \cap \Qzero$ and minor signs by Lemma~\ref{lem:invariance}. Finally, $\det(B) = (-1)^{m-1}\prod_{i=0}^{m-1} s_i = (-1)^{m-1}(-1)^{m-2t} = -1$.
\end{proof}

\section{Universal Jacobi Support Restriction, Dual-Frame Peeling, and Chain-Breaking}

Throughout the remainder of the paper, let $M = \begin{pmatrix} B & C \\ R & D \end{pmatrix} \in \Ezerof \cap \Qzero$ be in the canonical form \eqref{eq:blockM} of Theorem~\ref{thm:canonical}, with $\alpha = \{0, 1, \dots, m-1\}$ (indexed in $\Z_m$), $\bar{\alpha} = \{m, \dots, n-1\}$, $B_{i, i+1} = s_i \in \{+1, -1\}$, $C = (c_{i,k})$, $R = (r_{k,i})$, and $D = (d_{k,\ell})$. Denote $I_+ = \{i \in \Z_m : s_i = +1\}$ and $I_- = \{i \in \Z_m : s_i = -1\}$.

Since $(B^{-1})_{i+1, i} = s_i$, the PPT $M^\alpha = \begin{pmatrix} B^{-1} & \widetilde{C} \\ \widetilde{R} & \widetilde{D} \end{pmatrix}$ has entries
\begin{equation}\label{eq:Ma_entries}
  \tilde{c}_{i+1, k} = -s_i c_{i, k}, \qquad
  \tilde{r}_{k, i} = s_i r_{k, i+1}, \qquad
  \widetilde{D} = D - R B^{-1} C,
\end{equation}
for all $i \in \Z_m$ and $k \in \bar{\alpha}$. For $m \ge 3$, the $2 \times 2$ principal minors of $M$ and $M^\alpha$ give
\begin{align}
  \det M_{\{i, k\}} &= -c_{i, k} r_{k, i} \ge 0, \label{eq:2x2_M} \\
  \det(M^\alpha)_{\{i-1, k\}} &= -\tilde{c}_{i-1, k} \tilde{r}_{k, i-1} = s_{i-2} s_{i-1} c_{i-2, k} r_{k, i} \ge 0. \label{eq:2x2_Ma}
\end{align}

\begin{proof}[Proof of Theorem~\ref{thm:jacobi}]
Fix $k \in \bar{\alpha}$ and $i, u \in \Z_m$ with $u \notin \{i, i-1, i-2\} \pmod{m}$. Let $\gamma = \{i, i+1, \dots, u\} \subseteq \alpha$ be the contiguous cyclic arc from $i$ to $u$. Then $2 \le |\gamma| \le m - 2$, so both $|\gamma \cup \{k\}| = |\gamma| + 1 \le m - 1$ and $|(\alpha \setminus \gamma) \cup \{k\}| = (m - |\gamma|) + 1 \le m - 1$. By Theorem~\ref{thm:canonical}, $\det M_{\gamma\cup\{k\}} \ge 0$ and $\det(M^\alpha)_{(\alpha\setminus\gamma)\cup\{k\}} \ge 0$. By the standard PPT minor identity (see \cite{CPS2009}),
\[
  \det(M^\alpha)_{(\alpha\setminus\gamma)\cup\{k\}} = \frac{\det M_{\gamma\cup\{k\}}}{\det B} = -\det M_{\gamma\cup\{k\}},
\]
forcing $\det M_{\gamma\cup\{k\}} = 0$. Expanding $M_{\gamma\cup\{k\}}$ backward from row $u$ (whose sole nonzero entry is $c_{u,k}$) yields
\[
  0 = \det M_{\gamma\cup\{k\}} = (-1)^{|\gamma|}\Bigl(\prod_{j \in \gamma\setminus\{u\}} s_j\Bigr) c_{u,k} r_{k,i},
\]
so $c_{u,k} r_{k,i} = 0$. For any $i \in \Z_m$, the allowed set for $u$ is the $3$-arc $\{i-2, i-1, i\}$. Thus: if $|N_k| = 1$, $|\supp(c_{\bullet,k})| \le 3$; if $|N_k| = 2$, the intersection of two distinct $3$-arcs on $\Z_m$ ($m \ge 4$) has size at most $2$ (specifically, size $2$ when the two indices are adjacent, or when $m = 4$ and the two indices are opposite $\{i, i+2\}$ with intersection $\{i, i+2\}$; size $1$ when $m \ge 5$ and the indices are at cyclic distance $2$; and size $0$ otherwise); if $|N_k| = 3$, the intersection of three distinct $3$-arcs has size at most $1$ (size $1$ iff three consecutive indices, size $0$ otherwise); and if $|N_k| \ge 4$, the intersection is empty. Hence $|N_k| + |\supp(c_{\bullet,k})| \le 4$ whenever $N_k \neq \emptyset$ and $c_{\bullet,k} \neq 0$, and $|\supp(r_{k,\bullet})| \ge 4 \implies c_{\bullet,k} = 0$. The symmetric statement follows identically.
\end{proof}

For each $j \in I_+$ and real $t > 0$, define the test vector $q^{(j)}(t) \in \R^n$ by
\begin{equation}\label{eq:test_q}
  q^{(j)}_j(t) = -1, \qquad
  q^{(j)}_\ell(t) = +1 \;\; (\ell \in \alpha \setminus \{j\}), \qquad
  q^{(j)}_k(t) = t \;\; (k \in \bar{\alpha}).
\end{equation}
Partition $\alpha = \Z_m$ into $k_+$ pairwise disjoint chains ending at $I_+$: for each $j \in I_+$, let $p(j)$ and $\mathrm{next}(j)$ denote the preceding and following elements of $I_+$ along the cycle, and set
\[
  P(j) := \{p(j) + 1, \, p(j) + 2, \, \dots, \, j\} \subseteq \Z_m.
\]

\begin{lemma}[Feasibility, Dual Support, Peeling, and Chain-Breaking]\label{lem:toolkit}
For every $j \in I_+$ and every sufficiently large $t > 0$:
\begin{enumerate}
  \item $q^{(j)}(t)$ is feasible for $\LCP(q^{(j)}(t), M)$, so $\LCP(q^{(j)}(t), M)$ admits a solution $(z, w)$, and every solution satisfies $z_{\bar{\alpha}} \neq 0$, $z_j > 0$, and $w_j = 0$.
  \item $(q^{(j)}(t))^\alpha_{\ell+1} = -s_\ell q^{(j)}_\ell(t)$ for $\ell \in \Z_m$, and $(q^{(j)}(t))^\alpha_k > 0$ for $k \in \bar{\alpha}$. If $k_+ = 2$ with $I_+ = \{a, b\}$ and $j = b$, then $(q^{(b)}(t))^\alpha$ has its unique negative entry at $a + 1$, forcing $w_{a+1} > 0$ and $z_{a+1} = 0$.
  \item \emph{(Chain-Peeling)} There exist $k \in \bar{\alpha} \cap \supp(z)$ and $u \in P(j)$ with $r_{k, u} < 0$.
  \item \emph{(Chain-Breaking)} There exists $k \in \bar{\alpha}$ such that either $c_{j, k} > 0$ or $c_{u, k} < 0$ for some $u \in P(\mathrm{next}(j))$. Consequently, since each nonzero entry $c_{u,k}$ of $C$ can witness at most one index $j \in I_+$ (a positive entry $c_{u,k} > 0$ witnesses at most $j = u \in I_+$, while a negative entry $c_{u,k} < 0$ with $u \in P(j')$ witnesses at most $j = p(j')$), one has
  \begin{equation}\label{eq:C_count}
    \sum_{k \in \bar{\alpha}} |\supp(c_{\bullet, k})| \ge k_+.
  \end{equation}
  \item \emph{(Boundedness Barrier)} If for every sufficiently large $t > 0$ there is a solution $(z, w)$ of $\LCP(q^{(j)}(t), M)$ whose $R z_\alpha$ is bounded independently of $t$ (in particular, if $z_\alpha$ is bounded independently of $t$), a contradiction arises with $D \in \Ezero$.
\end{enumerate}
\end{lemma}

\begin{proof}
(1) Take $z = 2e_{j+1}$; then $Mz + q^{(j)}(t) > 0$ for $t > 2\|R\|_\infty$, so $M \in \Qzero$ yields a solution $(z, w)$. If $z_{\bar{\alpha}} = 0$, then $w_j = z_{j+1} - 1 \ge 0 \implies z_{j+1} \ge 1 > 0$; propagating along $P(\mathrm{next}(j))$ via $w_\ell = -z_{\ell+1} + 1 = 0 \implies z_{\ell+1} = 1 > 0$ reaches $j' = \mathrm{next}(j) \in I_+ \setminus \{j\}$ with $z_{j'} > 0$, contradicting $w_{j'} = z_{j'+1} + 1 \ge 1 > 0$. Thus $z_{\bar{\alpha}} \neq 0$. If $z_j = 0$, then $(Mz)_\ell = -q^{(j)}_\ell(t) < 0$ for all $\ell \in \supp(z)$, contradicting $M \in \Ezero$; hence $z_j > 0$ and $w_j = 0$.

(2) Immediate from $(B^{-1})_{\ell+1,\ell} = s_\ell$. When $I_+ = \{a, b\}$ and $j = b$, $-s_\ell q^{(b)}_\ell = +1$ for all $\ell \neq a$ and $-s_a q^{(b)}_a = -1$, so applying the $\Ezero$ argument of (1) to the pivoted solution $(u, v)$ of $\LCP((q^{(b)}(t))^\alpha, M^\alpha)$ gives $u_{a+1} = w_{a+1} > 0$ and $z_{a+1} = 0$.

(3) Set $z' = z - \sum_{u \in P(j)} z_u e_u \ge 0$. Since $z_{\bar{\alpha}} \neq 0$, $z' \neq 0$. For every $\ell \in \supp(z') \cap \alpha$, $q^{(j)}_\ell(t) = +1$ and $B_{\ell, u} \ge 0$ for all $u \in P(j)$ (the sole nonzero entry of $B$ from $P(j)$ to $\alpha \setminus P(j)$ is $B_{p(j), p(j)+1} = +1$), so $(Mz')_\ell = -1 - \sum_{u \in P(j)} z_u B_{\ell,u} \le -1 < 0$. By $M \in \Ezero$, there exists $k \in \supp(z') \cap \bar{\alpha}$ with $(Mz')_k = -t - \sum_{u \in P(j)} z_u r_{k,u} \ge 0$, forcing $r_{k,u} < 0$ for some $u \in P(j)$.

(4) Suppose for contradiction that $c_{j,k} \le 0$ and $c_{u,k} \ge 0$ for all $k \in \bar{\alpha}$ and all $u \in P(\mathrm{next}(j))$. Then in any solution $(z, w)$ of $\LCP(q^{(j)}(t), M)$, $w_j = z_{j+1} + \sum_k c_{j,k} z_k - 1 \ge 0$ forces $z_{j+1} \ge 1 > 0$. By induction along $u \in P(\mathrm{next}(j)) \setminus \{\mathrm{next}(j)\}$, $z_u > 0 \implies w_u = 0 \implies z_{u+1} = 1 + \sum_k c_{u,k} z_k \ge 1 > 0$. At $j' = \mathrm{next}(j)$, we obtain $z_{j'} > 0$ while $w_{j'} = z_{j'+1} + \sum_k c_{j',k} z_k + 1 \ge 1 > 0$, contradicting $z_{j'} w_{j'} = 0$. The counting bound \eqref{eq:C_count} follows immediately because the chains $P(j')$ ($j' \in I_+$) partition $\alpha$ and each entry $c_{u,k}$ has a single sign.

(5) Since $z_{\bar{\alpha}} \neq 0$, for every $k \in \supp(z_{\bar{\alpha}})$ complementarity gives $w_k = 0$, so $(D z_{\bar{\alpha}})_k = -(R z_\alpha)_k - t < 0$ once $t$ exceeds the uniform bound on $\|R z_\alpha\|_\infty$, contradicting $D \in \Ezero$.
\end{proof}

\begin{proof}[Proof of Theorem~\ref{thm:adj}]
If $k_+ = 2$ and $I_+ = \{a, b\}$ with $b \equiv a + 1 \pmod{m}$, then Lemma~\ref{lem:toolkit}(1)--(2) on $q^{(b)}(t)$ gives both $w_b = 0$ and $w_{a+1} = w_b > 0$, a contradiction. For $m = 1$, $B = (-1) \notin \Ezero$; for $m \in \{2, 3\}$, $k_+ = 2$ forces the two positive edges to be adjacent modulo $m$. Thus $m(A) \ge 4$.
\end{proof}

We also record the Cut-Peeling Lemma for $k_+ = 2$, which works for arbitrary $n \ge m + 1$.

\begin{lemma}[Cut-Peeling for $k_+ = 2$]\label{lem:cut_peeling}
Suppose $k_+ = 2$ with $I_+ = \{a, b\}$. Write $P(b) = \{v_0, v_1, \dots, v_\ell\}$ with $v_0 = a + 1$, $v_\ell = b$, and $\ell \ge 1$. Then there exist $k_1, k_2 \in \bar{\alpha}$ and $u_1, u_2 \in P(b) \setminus \{a + 1\}$ such that $r_{k_1, u_1} < 0$ and $r_{k_2, u_2} > 0$. In particular, if $\ell = 1$ (so $b = a + 2$), then $u_1 = u_2 = b$ and $k_1 \neq k_2$.
\end{lemma}

\begin{proof}
Let $(z, w)$ solve $\LCP(q^{(b)}(t), M)$ for large $t > 0$, and let $(u, v)$ be the pivoted solution to $\LCP((q^{(b)}(t))^\alpha, M^\alpha)$. By Lemma~\ref{lem:toolkit}(1)--(2), $z_{v_0} = 0$, $w_{v_0} > 0$, $z_{v_\ell} > 0$, and $w_{v_\ell} = 0$. Let $r \in \{0, \dots, \ell - 1\}$ be the smallest index with $z_{v_{r+1}} > 0$, so $z_{v_0} = \dots = z_{v_r} = 0$ and $w_{v_{r+1}} = 0$.
\begin{itemize}
  \item In $M$, set $z' = z - \sum_{i=r+1}^\ell z_{v_i} e_{v_i} \ge 0$. Since $z_{\bar{\alpha}} \neq 0$, $z' \neq 0$. For every $j \in \supp(z') \cap \alpha$, $j \notin \{v_r, \dots, v_\ell\}$ (as $z_{v_r} = 0$), so $B_{j, v_i} = 0$ for all $i \in \{r+1, \dots, \ell\}$ and $q^{(b)}_j = +1$; hence $(Mz')_j = -1 < 0$. By $M \in \Ezero$, there exist $k_1 \in \bar{\alpha}$ and $u_1 \in \{v_{r+1}, \dots, v_\ell\} \subseteq P(b) \setminus \{v_0\}$ with $r_{k_1, u_1} < 0$.
  \item In $M^\alpha$, set $u' = u - \sum_{i=0}^r u_{v_i} e_{v_i} \ge 0$. Since $u_{\bar{\alpha}} = z_{\bar{\alpha}} \neq 0$, $u' \neq 0$. For every $j \in \supp(u') \cap \alpha$, $j \notin \{v_0, \dots, v_{r+1}\}$ (as $u_{v_{r+1}} = w_{v_{r+1}} = 0$), so $(B^{-1})_{j, v_i} = 0$ for all $i \in \{0, \dots, r\}$ and $(q^{(b)}(t))^\alpha_j = +1$; hence $(M^\alpha u')_j = -1 < 0$. By $M^\alpha \in \Ezero$, there exist $k_2 \in \bar{\alpha}$ and $v_i \in \{v_0, \dots, v_r\}$ with $\tilde{r}_{k_2, v_i} < 0$. Since $i \le r \le \ell - 1$ implies $s_{v_i} = -1$, \eqref{eq:Ma_entries} gives $\tilde{r}_{k_2, v_i} = -r_{k_2, v_{i+1}} < 0$, so $r_{k_2, u_2} > 0$ for $u_2 = v_{i+1} \in P(b) \setminus \{v_0\}$.
\end{itemize}
If $\ell = 1$, then $P(b) \setminus \{v_0\} = \{b\}$, forcing $u_1 = u_2 = b$ and hence $k_1 \neq k_2$.
\end{proof}

\section{Elimination of \texorpdfstring{$n = m(A) + 1$}{n = m(A) + 1} and \texorpdfstring{$n = m(A) + 2$}{n = m(A) + 2}}

\begin{proof}[Proof of Theorem~\ref{thm:n_eq_m1}]
Suppose $A \in \Ezerof \cap \Qzero \setminus \Pzero$ has order $n$ with $m = m(A) = n - 1$ (while $m = n$ is excluded by Lemma~\ref{lem:toolkit}(1)). Pass to the canonical form \eqref{eq:blockM} with $m \ge 4$ and $\bar{\alpha} = \{m\}$, writing $C = c \in \R^m$ and $R = r^T \in \R^{1\times m}$. Since $c \neq 0$ (by Lemma~\ref{lem:toolkit}(4)), Theorem~\ref{thm:jacobi} implies $|\supp(r)| \le 3$.
\begin{itemize}
  \item If $k_+ \ge 4$, Lemma~\ref{lem:toolkit}(3) produces at least one $r_u < 0$ in each of the $k_+$ disjoint chains $P(j)$, giving $|\supp(r)| \ge k_+ \ge 4$, a contradiction.
  \item If $k_+ = 2$ with $I_+ = \{a, b\}$, Lemma~\ref{lem:cut_peeling} applied to $P(b)$ and $P(a)$ (with $|\bar{\alpha}| = 1$) forces $|P(b)| \ge 3$, $|P(a)| \ge 3$, and at least two distinct entries of $\supp(r)$ (one negative, one positive) in each of the disjoint sets $P(b) \setminus \{a+1\}$ and $P(a) \setminus \{b+1\}$, giving $|\supp(r)| \ge 4$, a contradiction. \qedhere
\end{itemize}
\end{proof}

We now prove our main all-$n$ strengthening, Theorem~\ref{thm:n_eq_m2}.

\begin{proof}[Proof of Theorem~\ref{thm:n_eq_m2}]
By Theorems~\ref{thm:adj} and~\ref{thm:n_eq_m1}, it suffices to rule out $n = m + 2$ for every $m \ge 4$. Pass to the canonical form \eqref{eq:blockM} with $\bar{\alpha} = \{m, m+1\}$. For $k \in \{m, m+1\}$, denote $S_k := \supp(r_{k,\bullet})$ and $N_k := \{i \in \Z_m : r_{k,i} < 0\}$. We split into $k_+ \ge 4$ and $k_+ = 2$.

\medskip
\noindent\textbf{Part 1: $k_+ \ge 4$.}
By Lemma~\ref{lem:toolkit}(3), $N_m \cup N_{m+1}$ intersects all $k_+$ disjoint chains $P(j)$ ($j \in I_+$), so $|N_m| + |N_{m+1}| \ge k_+$.

\emph{Subcase 1a: $c_{\bullet, m} \neq 0$ and $c_{\bullet, m+1} \neq 0$.}
First note that $N_m \neq \emptyset$ and $N_{m+1} \neq \emptyset$: if, say, $N_{m+1} = \emptyset$, then $|N_m| \ge k_+ \ge 4$, forcing $c_{\bullet, m} = 0$ by Theorem~\ref{thm:jacobi}, a contradiction. By Theorem~\ref{thm:jacobi}, $|N_k| + |\supp(c_{\bullet, k})| \le 4$ for each $k \in \{m, m+1\}$. Summing over $k \in \{m, m+1\}$ and applying \eqref{eq:C_count} yields
\[
  2 k_+ \le \bigl(|N_m| + |N_{m+1}|\bigr) + \bigl(|\supp(c_{\bullet, m})| + |\supp(c_{\bullet, m+1})|\bigr) \le 4 + 4 = 8.
\]
This immediately rules out all $k_+ \ge 6$, and for $k_+ = 4$ forces equality throughout: $N_m \cap N_{m+1} = \emptyset$, each element of $N_m \cup N_{m+1}$ lies in a distinct chain $P(j)$, $|N_k| + |\supp(c_{\bullet, k})| = 4$ for both $k \in \{m, m+1\}$, and every nonzero entry of $C$ witnesses a distinct index of $I_+$ in Lemma~\ref{lem:toolkit}(4). Up to swapping $m \leftrightarrow m+1$, $(|N_m|, |N_{m+1}|) \in \{(3, 1), (2, 2)\}$:
\begin{itemize}
  \item If $(|N_m|, |N_{m+1}|) = (3, 1)$, then $|\supp(c_{\bullet, m})| = 1$ requires $N_m = \{i, i+1, i+2\}$ (three consecutive indices) and $\supp(c_{\bullet, m}) = \{i\}$. Because $i, i+1, i+2$ lie in three distinct chains, $s_i = s_{i+1} = +1$. Then \eqref{eq:2x2_M} at $i$ gives $c_{i,m} \ge 0$, while \eqref{eq:2x2_Ma} at $i+2$ gives $s_i s_{i+1} c_{i,m} r_{m, i+2} \ge 0 \implies c_{i,m} \le 0$, forcing $c_{i,m} = 0$, a contradiction.
  \item If $(|N_m|, |N_{m+1}|) = (2, 2)$, then $|\supp(c_{\bullet, k})| = 2$ for each $k \in \{m, m+1\}$. Note first that the opposite-index case $m = 4$ with $N_k = \{i, i+2\}$ is impossible: for $m = 4$ and $k_+ = 4$ we have $s_0 = s_1 = s_2 = s_3 = +1$, and Theorem~\ref{thm:jacobi} restricts $\supp(c_{\bullet, k}) \subseteq \{i, i+2\}$, while \eqref{eq:2x2_M} and \eqref{eq:2x2_Ma} applied to $r_{k, i} < 0$ and $r_{k, i+2} < 0$ force both $c_{i, k} \ge 0, c_{i, k} \le 0$ and $c_{i+2, k} \ge 0, c_{i+2, k} \le 0$, giving $c_{\bullet, k} = 0$, a contradiction. 

  Therefore, by the arc-intersection characterization in the proof of Theorem~\ref{thm:jacobi}, each $N_k$ must be an adjacent pair $N_k = \{i_k, i_k+1\}$ with $\supp(c_{\bullet, k}) = \{i_k - 1, i_k\}$. Because $i_k, i_k+1$ lie in distinct chains, $s_{i_k} = +1$. By \eqref{eq:2x2_M} at $i_k$, $c_{i_k, k} > 0$ (witnessing $i_k \in I_+$). By \eqref{eq:2x2_Ma} at $i_k+1$, $s_{i_k-1} s_{i_k} c_{i_k-1, k} \le 0$. If $s_{i_k-1} = -1$, then $c_{i_k-1, k} > 0$ with $i_k - 1 \notin I_+$, which witnesses no index of $I_+$, contradicting equality in \eqref{eq:C_count}. Thus $s_{i_k-1} = +1$ and $c_{i_k-1, k} < 0$ (witnessing $p(i_k - 1)$). Now consider $j = i_m - 1 \in I_+$: since $P(\mathrm{next}(j)) = P(i_m) = \{i_m\}$, column $m$ does not witness $j$ (as $c_{i_m-1, m} < 0$ and $c_{i_m, m} > 0$). Hence column $m+1$ must witness $j$, requiring either $c_{i_m-1, m+1} > 0$ (so $i_{m+1} = i_m - 1 \implies i_{m+1}+1 = i_m \in N_m \cap N_{m+1}$) or $c_{i_m, m+1} < 0$ (so $i_{m+1} - 1 = i_m \implies i_{m+1} = i_m + 1 \in N_m \cap N_{m+1}$), contradicting $N_m \cap N_{m+1} = \emptyset$.
\end{itemize}

\emph{Subcase 1b: At least one column of $C$ is zero (without loss of generality $c_{\bullet, m} = 0$).}
By \eqref{eq:C_count}, $|\supp(c_{\bullet, m+1})| \ge k_+ \ge 4$, so Theorem~\ref{thm:jacobi} forces $r_{m+1, \bullet} = 0$ identically. Consequently $d_{m+1, m} < 0$: otherwise row $m+1$ of $M$ is entrywise nonnegative, so in $\LCP(q^{(j)}(t), M)$ we get $w_{m+1} \ge t > 0 \implies z_{m+1} = 0$, decoupling $w_\alpha = B z_\alpha + q^{(j)}_\alpha$ and contradicting Lemma~\ref{lem:toolkit}(1). Moreover, Lemma~\ref{lem:toolkit}(3) now places all $k_+$ chain-peeling negative entries in row $m$, so $|N_m| \ge k_+ \ge 4$.
\begin{itemize}
  \item \emph{Case $m \ge 6$:} For any $i, u \in \Z_m$ such that the contiguous cyclic arc $\gamma = \{i, i+1, \dots, u\}$ satisfies $3 \le |\gamma| \le m - 3$, both $|\gamma \cup \{m, m+1\}| \le m - 1$ and $|(\alpha \setminus \gamma) \cup \{m, m+1\}| \le m - 1$. By Theorem~\ref{thm:canonical} and the PPT minor identity,
  \[
    \det(M^\alpha)_{(\alpha\setminus\gamma)\cup\{m,m+1\}} = -\det M_{\gamma\cup\{m,m+1\}} \implies \det M_{\gamma\cup\{m,m+1\}} = 0.
  \]
  Because $c_{\bullet, m} = 0$ and $r_{m+1, \bullet} = 0$, expanding $M_{\gamma\cup\{m,m+1\}}$ backward from row $u$ (whose sole nonzero entry is $c_{u, m+1}$), then row $m+1$ (whose sole remaining nonzero entry is $d_{m+1, m}$), then rows $u-1, \dots, i$ yields the single monomial $\pm c_{u, m+1} r_{m, i} d_{m+1, m}$, which vanishes and hence gives $c_{u, m+1} r_{m, i} = 0$. Thus for every $u \in \supp(c_{\bullet, m+1})$, the set $N_m$ must lie inside the length-$5$ cyclic arc $W(u) := \{u - 1, u, u + 1, u + 2, u + 3\} \pmod{m}$. On $\Z_m$ with $m \ge 6$, the intersection of any three distinct length-$5$ arcs has cardinality at most $3$ (for $m = 6$, each complement is a distinct singleton, leaving $6 - 3 = 3$; for $m \ge 7$, only $W(u \pm 1)$ intersect $W(u)$ in $4$ elements, and $W(u-1) \cap W(u) \cap W(u+1)$ has size $3$). Since $|N_m| \ge 4$, we obtain $|\supp(c_{\bullet, m+1})| \le 2$, contradicting $|\supp(c_{\bullet, m+1})| \ge k_+ \ge 4$.
  \item \emph{Case $m \in \{4, 5\}$:} Here $k_+ = 4$, and up to cyclic shift either $m = 4$ with $(s_0, s_1, s_2, s_3) = (+1, +1, +1, +1)$ or $m = 5$ with $(s_0, s_1, s_2, s_3, s_4) = (-1, +1, +1, +1, +1)$. Direct expansion of the $3 \times 3$ minors of $M$ and $M^\alpha$ (using $c_{\bullet,m} = 0$, $r_{m+1,\bullet} = 0$, and $d_{m+1,m} < 0$) gives
  \begin{align*}
    \det M_{\{i, m, m+1\}} &= c_{i, m+1} r_{m, i} d_{m+1, m} \ge 0 \implies c_{i, m+1} r_{m, i} \le 0, \\
    \det(M^\alpha)_{\{i+1, m, m+1\}} &= -s_i s_{i+1} c_{i, m+1} r_{m, i+2} d_{m+1, m} \ge 0 \implies s_i s_{i+1} c_{i, m+1} r_{m, i+2} \ge 0.
  \end{align*}
  In both cases, $r_{m, 2} < 0$, $r_{m, 3} < 0$, and $r_{m, 4 \bmod m} < 0$ (for $m = 4$, all four $r_{m, i} < 0$; for $m = 5$, $P(2)=\{2\}, P(3)=\{3\}, P(4)=\{4\}$). Since $s_2 = s_3 = +1$, the two inequalities at $i = 2$ give $c_{2, m+1} \ge 0$ and $c_{2, m+1} \le 0$, so $c_{2, m+1} = 0$; while the first inequality at $i = 3$ gives $c_{3, m+1} \ge 0$. Because $P(\mathrm{next}(2)) = P(3) = \{3\}$, this violates Lemma~\ref{lem:toolkit}(4) at $j = 2$.
\end{itemize}

\medskip
\noindent\textbf{Part 2: $k_+ = 2$.}
Write $I_+ = \{a, b\}$ with chains $P(b) = \{a+1, \dots, b\}$ and $P(a) = \{b+1, \dots, a\}$, each of length at least $2$ by Theorem~\ref{thm:adj}.

\emph{Subcase 2a: At least one column of $C$ is zero (without loss of generality $c_{\bullet, m} = 0$).}
Consider any solution $(z, w)$ of $\LCP(q^{(b)}(t), M)$ (for large $t > 0$). By Lemma~\ref{lem:toolkit}(1)--(2), $z_{a+1} = 0$ and $w_b = 0$.
\begin{itemize}
  \item If $c_{u, m+1} < 0$ for some $u \in I_-$, then $w_u = -z_{u+1} + c_{u, m+1} z_{m+1} + 1 \ge 0$ forces $z_{m+1} \le 1/|c_{u, m+1}|$.
  \item Otherwise $c_{u, m+1} \ge 0$ for all $u \in I_-$. Then Lemma~\ref{lem:toolkit}(4) at $j = b$ requires either $c_{b, m+1} > 0$ (in which case $w_b = z_{b+1} + c_{b, m+1} z_{m+1} - 1 = 0$ forces $z_{m+1} \le 1/c_{b, m+1}$) or $c_{a, m+1} < 0$ (in which case $w_a = z_{a+1} + c_{a, m+1} z_{m+1} + 1 = c_{a, m+1} z_{m+1} + 1 \ge 0$ forces $z_{m+1} \le 1/|c_{a, m+1}|$).
\end{itemize}
In every case $z_{m+1}$ is bounded independently of $t$. It follows that $z_{u+1} \le 1 + c_{u, m+1} z_{m+1} = O(1)$ for every $u \in I_-$ (from $w_u \ge 0$), $z_{a+1} = 0$, and $z_{b+1} = 1 - c_{b, m+1} z_{m+1} = O(1)$ (from $w_b = 0$). Thus $z_\alpha$ is bounded independently of $t$, contradicting Lemma~\ref{lem:toolkit}(5).

\emph{Subcase 2b: $c_{\bullet, m} \neq 0$ and $c_{\bullet, m+1} \neq 0$, and at least one chain has length $2$ (without loss of generality $|P(b)| = 2$, so $b = a + 2$).}
By Lemma~\ref{lem:cut_peeling} on $P(b) = \{a+1, b\}$, after swapping $m \leftrightarrow m+1$ if necessary we have $r_{m, b} < 0$ and $r_{m+1, b} > 0$. By Theorem~\ref{thm:jacobi}, for both $k \in \{m, m+1\}$,
\[
  \supp(c_{\bullet, k}) \subseteq \{b - 2, b - 1, b\} = \{a, a + 1, b\}.
\]
Moreover, \eqref{eq:2x2_M} and \eqref{eq:2x2_Ma} at $i = b$ (where $s_a s_{a+1} = (+1)(-1) = -1$) yield
\begin{equation}\label{eq:len2_signs}
  c_{a, m} \ge 0, \quad c_{b, m} \ge 0, \qquad c_{a, m+1} \le 0, \quad c_{b, m+1} \le 0.
\end{equation}
If $m = 4$, then $|P(a)| = 2$ too ($a = b + 2$), so Lemma~\ref{lem:cut_peeling} on $P(a)$ gives $r_{m, a} r_{m+1, a} < 0$, forcing $c_{a+1, m} = c_{a+1, m+1} = 0$ by Theorem~\ref{thm:jacobi}; if $r_{m, a} > 0, r_{m+1, a} < 0$ then \eqref{eq:2x2_M}--\eqref{eq:2x2_Ma} reverse \eqref{eq:len2_signs} and force $C = 0$, while if $r_{m, a} < 0, r_{m+1, a} > 0$ we fall into the $c_{a+1, m} = c_{a+1, m+1} = 0$ case below.

If $m \ge 5$, then $|P(a)| \ge 3$ and $P(a) \setminus \{b+1\} = \{b+2, \dots, a\}$. Note that if $c_{a+1, k} \neq 0$, then Theorem~\ref{thm:jacobi} gives $S_k \subseteq \{a+1, b, b+1\}$, which is disjoint from $P(a) \setminus \{b+1\}$. By Lemma~\ref{lem:cut_peeling} on $P(a)$, $S_m \cup S_{m+1}$ must contain at least two entries (one negative, one positive) in $P(a) \setminus \{b+1\}$. Hence at least one of $c_{a+1, m}, c_{a+1, m+1}$ is zero. If exactly one is nonzero (say $c_{a+1, m} \neq 0, c_{a+1, m+1} = 0$), then $S_{m+1}$ must contain $b$ together with two distinct indices $v_1, v_2 \in \{b+2, \dots, a\}$; however, no $3$-arc on $\Z_m$ ($m \ge 5$) containing $b$ contains more than one index of $\{b+2, \dots, a\}$, forcing $c_{\bullet, m+1} = 0$ by Theorem~\ref{thm:jacobi}, a contradiction.

Thus in all cases $m \ge 4$ we have $c_{a+1, m} = c_{a+1, m+1} = 0$, and therefore $c_{u, k} = 0$ for every $u \in I_-$ and $k \in \{m, m+1\}$. Consequently, in both $\LCP(q^{(a)}(t), M)$ and $\LCP(q^{(b)}(t), M)$, every $u \in I_-$ satisfies $w_u = -z_{u+1} + 1 \ge 0$, so $z_{u+1} \in [0, 1]$. Now examine $c_{b, m+1} \le 0$ and $c_{a, m+1} \le 0$:
\begin{itemize}
  \item If $c_{b, m+1} = 0$, then in $\LCP(q^{(b)}(t), M)$ we have $z_{a+1} = 0$ and $w_b = z_{b+1} + c_{b, m} z_m - 1 = 0 \implies z_{b+1} \le 1$, so $z_\alpha \in [0, 1]^m$, contradicting Lemma~\ref{lem:toolkit}(5).
  \item Symmetrically, if $c_{a, m+1} = 0$, then in $\LCP(q^{(a)}(t), M)$ we have $z_{b+1} = 0$ and $w_a = z_{a+1} + c_{a, m} z_m - 1 = 0 \implies z_{a+1} \le 1$, contradicting Lemma~\ref{lem:toolkit}(5).
  \item If $c_{b, m+1} < 0$ and $c_{a, m+1} < 0$, then Theorem~\ref{thm:jacobi} gives $S_{m+1} \subseteq \{b, b+1, b+2\} \cap \{a, a+1, b\} = \{b\}$ (for $m \ge 5$) or $\{a, b\}$ (for $m = 4$); in either case $r_{m+1, b+1} = 0$. Now consider $\LCP(q^{(b)}(t), M)$: all coordinates of $z_\alpha$ except possibly $z_{b+1}$ lie in $[0, 1]$.
  \begin{itemize}
    \item If $c_{a, m} > 0$, then Theorem~\ref{thm:jacobi} gives $S_m \subseteq \{a, a+1, b\}$, so $r_{m, b+1} = 0$ too; hence $R z_\alpha$ is independent of $z_{b+1}$ and bounded independently of $t$, contradicting Lemma~\ref{lem:toolkit}(5).
    \item If $c_{a, m} = 0$, then $w_a = z_{a+1} + c_{a, m+1} z_{m+1} + 1 = c_{a, m+1} z_{m+1} + 1 \ge 0$ implies $z_{m+1} \le 1/|c_{a, m+1}|$, and $w_b = 0$ then gives $z_{b+1} = 1 - c_{b, m} z_m - c_{b, m+1} z_{m+1}$, whence $z_{b+1} \le 1 + |c_{b, m+1}|/|c_{a, m+1}|$. Thus $z_\alpha$ is bounded independently of $t$, contradicting Lemma~\ref{lem:toolkit}(5).
  \end{itemize}
\end{itemize}

\emph{Subcase 2c: $c_{\bullet, m} \neq 0$, $c_{\bullet, m+1} \neq 0$, and $|P(a)| \ge 3$, $|P(b)| \ge 3$ (hence $m \ge 6$).}
By Theorem~\ref{thm:jacobi}, each of $S_m$ and $S_{m+1}$ lies inside a contiguous $3$-arc of $\Z_m$. By Lemma~\ref{lem:cut_peeling}, $P(b) \setminus \{a+1\}$ contains both a negative and a positive entry of $R$, and $P(a) \setminus \{b+1\}$ likewise contains both a negative and a positive entry of $R$. Notice that $P(b) \setminus \{a+1\}$ and $P(a) \setminus \{b+1\}$ are separated on $\Z_m$ by the two gap vertices $a+1$ and $b+1$; hence a $3$-arc can intersect both sets only if it is $\{a, a+1, a+2\}$ or $\{b, b+1, b+2\}$.
\begin{itemize}
  \item \emph{Neither $S_m$ nor $S_{m+1}$ crosses a gap.} Then, after swapping $m \leftrightarrow m+1$ if necessary, $S_m \cap (P(b) \setminus \{a+1\})$ contains two distinct entries $u_1, u_2$ with $r_{m, u_1} < 0 < r_{m, u_2}$ (and $S_m \cap (P(a) \setminus \{b+1\}) = \emptyset$), while $S_{m+1} \cap (P(a) \setminus \{b+1\})$ contains two distinct entries $v_1, v_2$ with $r_{m+1, v_1} < 0 < r_{m+1, v_2}$. Since $S_m$ lies in a $3$-arc and contains $\{u_1, u_2\} \subseteq P(b) \setminus \{a+1\} = \{a+2, \dots, b\}$, Theorem~\ref{thm:jacobi} forces $\supp(c_{\bullet, m}) \subseteq \{u_1-2, u_1-1, u_1\} \cap \{u_2-2, u_2-1, u_2\} \subseteq P(b) \setminus \{b\} \subseteq I_-$ (specifically, a subset of $\{i\}$ if $\{u_1, u_2\} = \{i, i+2\}$, or of $\{i-1, i\}$ if $\{u_1, u_2\} = \{i, i+1\}$). In both cases, all edges in $P(b) \setminus \{b\}$ have sign $-1$, so \eqref{eq:2x2_M} and \eqref{eq:2x2_Ma} force all nonzero entries of $c_{\bullet, m}$ to have the same sign. Symmetrically, $\supp(c_{\bullet, m+1}) \subseteq P(a) \setminus \{a\} \subseteq I_-$ and all nonzero entries of $c_{\bullet, m+1}$ have the same sign. 

  By Lemma~\ref{lem:toolkit}(4) at $j = a$ (where $c_{a, k} = 0$ for both $k$), $c_{\bullet, m}$ must have a strictly negative entry $c_{u, m} < 0$ ($u \in P(b) \setminus \{b\} \subseteq I_-$), and hence $c_{\bullet, m} \le 0$; symmetrically $c_{\bullet, m+1} \le 0$ with some $c_{v, m+1} < 0$ ($v \in P(a) \setminus \{a\} \subseteq I_-$). Now in $\LCP(q^{(b)}(t), M)$: since $c_{u, m+1} = 0$ and $c_{v, m} = 0$, the inequalities $w_u = -z_{u+1} + c_{u,m} z_m + 1 \ge 0$ and $w_v = -z_{v+1} + c_{v,m+1} z_{m+1} + 1 \ge 0$ bound both $z_m \le 1/|c_{u,m}|$ and $z_{m+1} \le 1/|c_{v,m+1}|$. Hence $z_{j+1} \le 1$ for all $j \in I_-$, $z_{a+1} = 0$, and $w_b = z_{b+1} - 1 = 0 \implies z_{b+1} = 1$, contradicting Lemma~\ref{lem:toolkit}(5).
  \item \emph{At least one of $S_m, S_{m+1}$ (say $S_m$) crosses a gap, without loss of generality $\{a, a+2\} \subseteq S_m \subseteq \{a, a+1, a+2\}$.} By Theorem~\ref{thm:jacobi}, $\supp(c_{\bullet, m}) = \{a\}$. By \eqref{eq:2x2_M} at $a$ and \eqref{eq:2x2_Ma} at $a+2$ (where $s_a s_{a+1} = -1$), both $c_{a, m} r_{m, a} \le 0$ and $c_{a, m} r_{m, a+2} \le 0$, so $r_{m, a}$ and $r_{m, a+2}$ have the same sign (opposite to $c_{a, m}$). Consequently, to supply the opposite sign in both $P(b) \setminus \{a+1\}$ and $P(a) \setminus \{b+1\}$ (Lemma~\ref{lem:cut_peeling}), $S_{m+1}$ must also intersect both sets, hence must also cross a gap:
  \begin{itemize}
    \item If $S_{m+1}$ crosses at $a+1$, then $\supp(c_{\bullet, m+1}) = \{a\}$ as well, so every row in $\alpha \setminus \{a\}$ of $C$ is zero. In $\LCP(q^{(b)}(t), M)$, every $u \in I_-$ satisfies $w_u = -z_{u+1} + 1 \ge 0 \implies z_{u+1} \in [0, 1]$, $z_{a+1} = 0$, and $w_b = z_{b+1} - 1 = 0 \implies z_{b+1} = 1$, contradicting Lemma~\ref{lem:toolkit}(5).
    \item If $S_{m+1}$ crosses at $b+1$, then $\supp(c_{\bullet, m+1}) = \{b\}$, so $c_{u, k} = 0$ for all $u \in I_-$ and $k \in \{m, m+1\}$. If $c_{a, m} > 0$, then in $\LCP(q^{(a)}(t), M)$ we get $z_{b+1} = 0$, $z_{u+1} \in [0, 1]$ ($u \in I_-$), and $w_a = z_{a+1} + c_{a,m} z_m - 1 = 0 \implies z_{a+1} \le 1$, contradicting Lemma~\ref{lem:toolkit}(5). Symmetrically $c_{b, m+1} > 0$ is impossible via $\LCP(q^{(b)}(t), M)$. Thus $c_{a, m} < 0$ and $c_{b, m+1} < 0$. By \eqref{eq:2x2_M}--\eqref{eq:2x2_Ma}, this forces $r_{m, a}, r_{m, a+2} > 0$ and $r_{m+1, b}, r_{m+1, b+2} > 0$, leaving no negative entry of $R$ in $P(b) \setminus \{a+1\}$, contradicting Lemma~\ref{lem:cut_peeling}.
  \end{itemize}
\end{itemize}
All cases yield a contradiction, establishing $m(A) \le n - 3$.
\end{proof}

\begin{remark}[Counting Barrier at Codimension $p = n - m = 3$ and Symbolic Verification]\label{rem:boundary}
The proof of Theorem~\ref{thm:n_eq_m2} shows why the present support-counting argument does not extend directly to codimension $p = n - m = 3$ (first occurring at $n = 7, m = 4$). For $p = 2$, the bound $|N_k| + |\supp(c_{\bullet, k})| \le 4$ forces exact tightness $2k_+ \le 4p = 8$ at $k_+ = 4$, and two rows exhaust $\bar{\alpha}$ in the $k_+ = 2$ peeling step. At $p = 3$, the sum bound $2k_+ \le 4p = 12$ no longer constrains $k_+ = 4$: for instance, on the positive $4$-cycle at $n = 7$, the partition $N_4 = \{0\}$, $N_5 = \{1\}$, $N_6 = \{2, 3\}$ satisfies Theorem~\ref{thm:jacobi} while allowing strictly positive entries in every row of $C$, so the two-row sign-pigeonhole argument ceases to yield an immediate contradiction.

The accompanying ancillary script checks selected determinant, PPT, contiguous-arc Jacobi, and $3 \times 3$ minor identities, together with the arc-intersection bounds for $m \le 10$, in exact symbolic arithmetic.
\end{remark}


\begin{thebibliography}{10}

\bibitem{BM1997}
A.~K.~Biswas and G.~S.~R.~Murthy,
\emph{A chain condition for $Q_0$-matrices},
in \emph{Game Theoretical Applications to Economics and Operations Research}
(T.~Parthasarathy, B.~Dutta, J.~A.~M.~Potters, T.~E.~S.~Raghavan, D.~Ray, and A.~Sen, eds.),
Kluwer Academic Publishers, Dordrecht, 1997, pp.~149--152.

\bibitem{CPS2009}
R.~W.~Cottle, J.-S.~Pang, and R.~E.~Stone,
\emph{The Linear Complementarity Problem},
Classics in Applied Mathematics, Vol.~60, SIAM, Philadelphia, 2009 (originally Academic Press, 1992).

\bibitem{CS1983}
R.~W.~Cottle and R.~E.~Stone,
\emph{On the uniqueness of solutions to linear complementarity problems},
Math.\ Program.\ \textbf{27} (1983), 191--213.

\bibitem{Ghosh2026}
S.~Ghosh,
\emph{Sign propagation of principal minors in fully semimonotone matrices and Stone's conjecture},
arXiv:2608.22829v3, October 2026,
\url{https://arxiv.org/abs/2608.22829}.

\bibitem{MND2001}
S.~R.~Mohan, S.~K.~Neogy, and A.~K.~Das,
\emph{On the classes of fully copositive and fully semimonotone matrices},
Linear Algebra Appl.\ \textbf{323} (2001), 87--97.

\bibitem{MP1995}
G.~S.~R.~Murthy and T.~Parthasarathy,
\emph{Some properties of fully semimonotone, $Q_0$-matrices},
SIAM J.\ Matrix Anal.\ Appl.\ \textbf{16} (1995), no.~4, 1268--1286.

\bibitem{MP1998}
G.~S.~R.~Murthy and T.~Parthasarathy,
\emph{Fully copositive matrices},
Math.\ Program.\ \textbf{82} (1998), 401--411.

\bibitem{MPS2013}
G.~S.~R.~Murthy, T.~Parthasarathy, and R.~Sridhar,
\emph{On fully semimonotone matrices},
Int.\ Game Theory Rev.\ \textbf{15} (2013), no.~4, 1340036.

\bibitem{ND2005}
S.~K.~Neogy and A.~K.~Das,
\emph{Principal pivot transforms of some classes of matrices},
Linear Algebra Appl.\ \textbf{400} (2005), 243--252.

\bibitem{Stone1981}
R.~E.~Stone,
\emph{Geometric Aspects of the Linear Complementarity Problem},
Ph.D.\ Thesis, Department of Operations Research, Stanford University, 1981.

\end{thebibliography}
\end{document}